\documentclass{article}
\usepackage{amssymb,latexsym} 
\newcommand{\1}{{{\mathchoice {\rm 1\mskip-4mu l} {\rm 1\mskip-4mu l} 
{\rm 1\mskip-4.5mu l} {\rm 1\mskip-5mu l}}}} 

\newcommand{\R}{{\mathbb R}} 

\newcommand{\Q}{{\mathbb Q}}
\newcommand{\QH}{{\rm QH}} 
 
\newcommand{\SO}{{\rm SO}} 
\newcommand{\Flux}{{\rm Flux}} 
\newcommand{\Symp}{{\rm Symp}}

\newcommand{\GW}{{\rm GW}} 
\newcommand{\TJ}{{\widetilde{J}}}

\newcommand{\SL}{{\rm SL}} 
 
\newcommand{\tr}{{\rm tr}}

\newcommand{\PU}{{\rm PU}}

\newcommand{\ev}{{\rm ev}}

\newcommand{\TLl} {{\Tilde{\Ll}}}

\newcommand{\Spec}{{\rm Spec}} 
 
\newcommand{\Hh}{{\mathcal H}}

\newcommand{\Z}{{\mathbb Z}} 
\newcommand{\C}{{\mathbb C}}

\newcommand{\Xx}{{\mathcal {X}}}

\newcommand{\G}{{\rm G}} 
\newcommand{\p}{{\partial}} 
\newcommand{\al}{{\alpha}}

\newcommand{\Om}{{\Omega}} 
\newcommand{\om}{{\omega}}

\newcommand{\eps}{{\varepsilon}}

\newcommand{\ga}{{\gamma}} 
\newcommand{\Ga}{{\Gamma}} 
\newcommand{\ka}{{\kappa}} 
\newcommand{\la}{{\lambda}} 
\newcommand{\si}{{\sigma}}

\newcommand{\Aa}{{\mathcal A}}

\newcommand{\Dd}{{\mathcal D}} 
\newcommand{\Ll}{{\mathcal L}}

\newcommand{\Nn}{{\mathcal N}}

\newcommand{\Ss}{{\mathcal S}}

\newcommand{\La}{{\Lambda}}

\newcommand{\Si}{{\Sigma}} 
\newcommand{\Ham}{{\rm Ham}}

\newcommand{\IFF}{{\Longleftrightarrow}} 
 
\newcommand{\MS}{{\medskip}}

\newcommand{\NI}{{\noindent}}

\newcommand{\QED}{\hfill$\Box$\medskip}

\newcommand{\Tilde}{\widetilde}
\newcommand{\THam}{{\Tilde{\Ham}}}
\newcommand{\TSymp}{{\widetilde{\Symp}}}

\newcommand{\tphi}{{\Tilde{\phi}}}

\newcommand{\CP}{{\mathbb CP}} 

\newcommand{\Diff}{{\rm Diff}}

\newtheorem{theorem}{Theorem}[section] 
\newtheorem{cor}[theorem]{Corollary}

\newtheorem{lemma}[theorem]{Lemma} 
\newtheorem{question}[theorem]{Question} 
 
\newtheorem{prop}[theorem]{Proposition}

\begin{document} 

\title{A survey of the topological properties  of symplectomorphism groups}
\author{Dusa McDuff\thanks{Partially
supported by NSF grant DMS 0072512}
\\ State University of New York
at Stony Brook \\ (dusa@math.sunysb.edu)}

\date{January 28, 2003}

\maketitle
\MS\MS
\begin{center}  for Graeme Segal
    \end{center}
\MS\MS

\begin{abstract}
   The special structures that arise in symplectic topology (particularly 
Gromov--Witten invariants and  quantum homology) place as yet 
rather poorly 
understood restrictions on  the topological properties of 
symplectomorphism groups.  
This article surveys some recent work by Abreu, Lalonde, McDuff, Polterovich 
and Seidel,
concentrating particularly on the homotopy properties of the action of 
the group of  Hamiltonian symplectomorphisms  on 
the underlying manifold $M$.  It sketches the proof that the 
evaluation map $\pi_1(\Ham(M))\to \pi_1(M)$ given by 
$\{\phi_t\}\mapsto \{\phi_t(x_0)\}$ is trivial, as well as explaining 
similar vanishing results for the  action of the homology of $\Ham(M)$ 
on the homology of $M$.  Applications to Hamiltonian stability are 
discussed.
\end{abstract}

\section{Overview}

The special structures\footnotetext{
\NI
 keywords:  symplectomorphism group, 
 Hamiltonian group, symplectic bundles, quantum 
 homology, Seidel 
 representation\\
 \NI
 Mathematics Subject Classification 2000: 57R17, 53D35}
 that arise in symplectic topology (particularly 
Gromov--Witten invariants and  quantum homology) place as yet 
rather poorly 
understood restrictions on  the topological properties of 
symplectomorphism groups.  
This article surveys some recent work on this subject.
Throughout $(M,\om)$ will be a closed (ie compact and without boundary),
smooth
symplectic manifold of dimension $2n$, 
unless it is explicitly mentioned otherwise.  Background information 
and more references can be found in~\cite{MS,JHOL,Pbk}.

The {\bf symplectomorphism group} $\Symp(M,\om)$ 
consists of all 
diffeomorphisms $\phi: M\to M$ such that $\phi^*(\om) = \om$,
and is equipped with the $C^\infty$-topology, the topology of uniform 
convergence of all derivatives.  We will sometimes contrast this with 
the $C^0$ (i.e.  compact-open) topology.
The (path) connected component containing the identity  
is denoted $\Symp_0(M,\om)$. (Note that $\Symp$ is locally path connected.) 
This group $\Symp_0$
contains an important  normal subgroup called the {\bf Hamiltonian group} 
$\Ham(M,\om)$ whose elements are the 
time-$1$ maps of Hamiltonian flows. 
These are the flows $\phi_t^H, t\in [0,1],$ that at each time $t$ 
are tangent to
the symplectic gradient $X_t^H$ of the function $H_t:M\to \R$, i.e.
$$
\dot{\phi}_t^H = X_t^H, \qquad \om(X_t^H,\cdot) = -dH_t.
$$
When $H^1(M, \R) = 0$ the groups 
$\Ham$ and $\Symp_0$ coincide.
In general, there is a sequence of groups
and inclusions
$$
\Ham(M, \om)\hookrightarrow \Symp_0(M, \om)\hookrightarrow\Symp(M,\om)
\hookrightarrow\Diff^+(M),
$$
where $\Diff^+$ denotes the orientation preserving diffeomorphisms.
Our aim is to understand and contrast the properties of these groups.

We first give an overview of basic results on the group $\Symp_0$.
Then we describe results on the Hamiltonian group, showing how a 
vanishing theorem 
for its action on $H_*(M)$ implies various 
stability results.  Finally, we sketch the proof of this vanishing 
theorem.  It relies on properties of the Gromov--Witten invariants 
for sections of Hamiltonian fiber bundles over $S^2$, that can be summarized 
in the statement (essentially
due to Floer and Seidel) that there is a representation 
of $\pi_1(\Ham(M,\om))$ into  the automorphism group 
of the quantum homology ring of $M$.   The proof of the vanishing of
the evaluation 
map $\pi_1(\Ham(M))\to \pi_1(M)$ is easier: it relies on a 
\lq\lq stretching the neck'' argument, see Lemma~\ref{le:sect} below.
A different but also relatively easy proof of this fact
may be found in~\cite{JHOL}.

\subsection*{Basic facts}

We begin by listing some 
fundamentals.\MS

 \NI
 $\bullet$ {\it Dependence on the cohomology class of $\om$.}
 
 The groups $\Symp(M,\om)$ and $\Ham(M, \om)$ depend 
only on the diffeomorphism class of the form $\om$.  In particular, 
since Moser's argument implies that 
any path $\om_t, t\in [0,1],$ of {\it cohomologous} forms is 
 induced by an isotopy $\psi_t: M\to M$  of the underlying manifold (i.e. 
 $\psi_t^*(\om_t) = \om_0$, $\psi_0 = id$), the groups do not change their 
 topological or algebraic properties when $\om_t$ varies along such a 
 path.  However, as first notices by Gromov (see Proposition~\ref{prop:am2}
 below), changes in the cohomology class  $[\om]$ can cause 
 significant changes in the homotopy type of these groups.
\MS
 
\NI
 $\bullet$   {\it Stability properties of $\Symp(M)$ and 
 $\Symp_0(M)$.}
 
 By this we mean that if $G$ denotes either of these
 groups, there is a $C^1$-neighbourhood $\Nn(G)$ of $G$
 in $\Diff(M)$ that deformation retracts onto $G$.  This follows from 
 the Moser isotopy argument mentioned above. In the case $G = 
 \Symp(M)$, take
 $$
 \Nn(\Symp)= \{\phi\in \Diff(M): (1-t)\phi^*(\om) + t\om \mbox{ is 
 nondegenerate for } t\in [0,1]\}.
 $$
 By Moser, one can define for each such $\phi$ a unique  isotopy 
 $\psi_t$ (that 
 depends smoothly on $\phi^*(\om)$) such that 
 $\psi_t^*(t\phi^*(\om) + 
(1-t)\om) = \om$ for all $t$.  Hence $\phi\circ\psi_1\in \Symp(M)$.
 Similarly, when $G= \Symp_0(M)$ one can take $ \Nn(G)$ to be the
 identity component of $\Nn(\Symp)$.  Note also that these 
 neighborhoods are uniform with respect to $\om$.  For example, 
 given any compact subset $K$ of $\Symp_0(M, \om)$  there is a 
 $C^{\infty}$-neighbourhood $\Nn(\om)$ of $\om$ in the 
 space of all symplectic forms such that $K$ may be isotoped into 
 $\Symp_0(M, \om')$ for all $\om'\in \Nn(\om)$.  These 
 statements, that we sum up in the rubric {\bf symplectic stability},
 exhibit the flabbiness, or lack of local invariants, of symplectic 
 geometry.
 
  \MS

  The above two properties are \lq\lq soft'', i.e. they depend only on 
  the Moser argument.  By way of contrast, the next result is \lq\lq 
  hard'' and can be proved only by using some deep ideas, either from 
  variational calculus (Ekeland--Hofer), generating functions/wave 
  fronts (Eliashberg, Viterbo)   or $J$-holomorphic curves 
  (Gromov).\MS
  
 \NI
 $\bullet$ {\it The group $\Symp(M,\om)$ is $C^0$-closed in $\Diff(M)$.}
 
 This celebrated result of Eliashberg and Ekeland--Hofer 
 is known as {\bf symplectic rigidity} and
 is the basis 
 of  symplectic topology.  The proof shows that even though 
 one uses the first derivatives of $\phi$ in 
 saying that a diffeomorphism $\phi$ preserves $\om$, there is an
 invariant $c(U)$ (called a {\it symplectic capacity}) of an open 
 subset of a symplectic manifold that is
continuous with respect to the Hausdorff metric on sets and that is 
preserved by a diffeomorphism $\phi$ if and only if 
$\phi^*(\om) = \om$. (When $n$ is even,
one must slightly modify the previous statement to rule out 
the case $\phi^*(\om) = -\om$.)  There are 
several ways to define a suitable invariant $c$.   Perhaps the easiest
is to take Gromov's width:
$$
c(U) = \sup \{\pi r^2: B^{2n}(r)\mbox{ 
embeds symplectically in }U\}.
$$
Here  $B^{2n}(r)$ is the standard ball of radius $r$ in Euclidean 
space $\R^{2n}$ with the usual symplectic form $\om_0 = \sum_i 
dx_{2i-1}\wedge dx_{2i}$.\MS

It is unknown whether the identity component $\Symp_0(M)$ is 
$C^0$-closed in $\Diff(M)$.  In fact this may well not hold.  For example, 
it is quite possible that the group $\Symp^c(\R^{2n})$ of compactly supported 
symplectomorphisms of Euclidean space is disconnected when $n> 2$.
(When $n=2$ this group is contractible by Gromov~\cite{GRO}.)  Hence 
 for 
some closed manifold $M$ there might be an element in $\Symp(M)\setminus 
\Symp_0(M)$ that is supported 
in a Darboux neighbourhood $U$ (i.e. an open set symplectomorphic to an 
open ball in Euclidean space).  Such an element would be in the $C^0$-closure of 
$\Symp_0(M)$ since by conformal rescaling in $U$ one could isotop it to have 
support in an arbitrarily small neighbourhood of a point in $U$.\MS

We discuss related questions for the group $\Ham(M)$ in Section~2 below. 
Though less is known about the above questions, some very interesting 
new features appear.
Before doing that we shall give a brief summary of what is known about
the homotopy groups of $\Symp(M)$.

\subsection*{The homotopy type of $\Symp(M)$}

In dimension $2$ it follows from Moser's argument that 
$\Symp(M,\om)$ is homotopy equivalent to $\Diff^+$.  Thus 
$
\Symp(S^2)$ is homotopy equivalent to the rotation group $ \SO(3)$;
$\Symp_0(T^2)$ is homotopy equivalent to an 
extension of $\SL(2,\Z)$ by $T^2$; and for 
 higher genus the symplectomorphism group 
 is homotopy equivalent to the mapping class 
group.  In dimensions $4$ and above, almost nothing
is known about the homotopy type 
of $\Diff^+$.  On the other hand, there are some very special $4$-manifolds for 
which  the 
(rational) homotopy type of $\Symp$ is fully understood.
 The following 
results are due to Gromov~\cite{GRO}.  Here $\si_Y$  
denotes (the pullback to the product of) an area form on
the Riemann surface $Y$ with total area $1$.

\begin{prop}[Gromov]
\begin{description}\item[(i)]
 $\Symp^c(\R^4, \om_0)$ is contractible;
\item[(ii)]
$\Symp(S^2\times S^2, \si_{S^2} + \si_{S^2})$ is homotopy equivalent to the extension 
of $\SO(3)\times \SO(3)$ by $\Z/2\Z$ where this acts by interchanging the 
factors;
\item[(iii)]  $\Symp(\C P^2, \om_{\rm FS})$ is 
homotopy equivalent to $\PU(3)$, where $\om_{\rm FS}$ is the 
Fubini--Study K\"ahler form. 
\end{description}
\end{prop}

It is no coincidence that these results occur in dimension $4$.  The 
proofs use $J$-holomorphic spheres, and these give much more information in 
dimension $4$ because of positivity of intersections. 

In Abreu~\cite{Abr} and Abreu--McDuff~\cite{AM} 
these arguments are extended  to other symplectic forms 
and (some) other ruled surfaces.  
Here are the main results, stated for convenience for the product 
manifold $\Sigma\times S^2$ (though there are similar results for the 
nontrivial $S^2$ bundle over $\Si$.)  
Consider the following family\footnote
{
Using results of Taubes and Li--Liu, Lalonde--McDuff show in~\cite{LMc} that 
these are the {\it only} symplectic forms on $\Si\times S^2$ up to 
diffeomorphism.} of symplectic forms on $M_g = \Sigma_g\times S^2$
(where $g$ is  genus$(\Si)$):
$$
\om_\mu = \mu \si_\Si +  \si_{S^2},\qquad \mu > 0.
$$
Denote by $G_\mu^g$ the subgroup 
$$
 G_\mu^g:= \Symp(M_g,
\om_\mu)\cap \Diff_0(M_g)
$$
of the group of symplectomorphisms of $(M_g,
\om_\mu)$.  When $g > 0$, $\mu$ ranges over all positive numbers. 
However, when $g=0$ there is an extra symmetry ---  interchanging the two
spheres gives an isomorphism $G_\mu^0 \cong G_{1/\mu}^0$ --- and so we 
take $\mu \ge 1$.  Although it is not completely obvious, 
there is a natural homotopy class of maps  from  $G_\mu^g$ to $
G_{\mu + \eps}^g$ for all $\eps>0$.  To see this, let
$$
G_{[a,b]}^g = \bigcup_{\mu\in [a,b]} \; \{\mu\}\times G_\mu^g 
\;\;\subset \;\;\R\times \Diff(M_g).
$$
It is shown in~\cite{AM} that the inclusion $G_b^g\to G_{[a,b]}^g$
is a homotopy equivalence. Therefore we can take the map
 $G_\mu^g \to G_{\mu + \eps}^g$ to be the composite of the inclusion
$G_\mu^g\to G_{[\mu, \mu+\eps]}^g$ with a homotopy inverse 
$G_{[\mu, \mu+\eps]}^g\to G_{\mu+\eps}^g$. Another, more geometric
definition of this map is given in~\cite{Mcrs}.

\begin{prop}\label{prop:am1}  As $\mu\to \infty$,  the groups $G_\mu^g$
tend to a limit $G_\infty^g$ that
has the homotopy type  of the identity component $\Dd_0^g$ of the group 
of fiberwise diffeomorphisms
of $M_g = \Si_g\times S^2 \to \Si_g$.
\end{prop}

\begin{prop}\label{prop:am2}
When $\ell  < \mu \le \ell + 1$ for some integer
$\ell \ge 1$,  
$$
H^*(G_\mu^0, \Q) = \La(t, x, y)\otimes \Q[w_\ell], 
$$
where $\La(t,x,y)$ is an exterior algebra over $\Q$ with generators $t$ of
degree $1$, and $x,y$ of degree $3$ and $\Q[w_\ell]$ is the polynomial
algebra
on a generator $w_\ell$ of degree $4\ell$.
\end{prop}

In the above statement, the generators $x,y$ come from 
$H^*(G_1^0) = H^*(\SO(3)\times
\SO(3))$ and $t$ corresponds to an element in $\pi_1(G_\mu^0), \mu >
1$ found by Gromov in~\cite{GRO}.   
Thus the subalgebra $\La(t, x, y)$ is the pullback of $H^*(\Dd_0^0, \Q)$
under the map $G_\mu^0 \to \Dd_0^0$.    The other generator
$w_\ell$ is fragile, in the sense that the corresponding element in 
homology disappears (i.e. becomes null
homologous) when $\mu$ increases.   It is dual to an element in 
$\pi_{4\ell}$ that is a higher order Samelson product and hence gives 
rise to a relation  (rather than a new generator) in the cohomology of 
the classifying space.  Indeed, when $\ell  < \mu \le \ell + 1$,
$$
H^*(BG^0_\mu) \cong \frac{\Q[T,X,Y]}{\{T(X-Y)\dots(\ell^2 X-Y)= 0\}},
$$
where the classes $T,X,Y$ have dimensions $2,4,4$ respectively and 
are the deloopings of $t,x,y$.

Anjos~\cite{A}  calculated the 
full homotopy type of $G_\mu^0$ for $1< \mu \le 2$. Her results
has been sharpened in Anjos--Granja~\cite{AG} where it is shown that
this group has the homotopy type of the pushout of the following 
diagram in the category of topological groups:
$$
\begin{array}{ccc} \SO(3)&\stackrel{diag}\longrightarrow & \SO(3)\times \SO(3)\\
    \downarrow& &\\
    S^1\times \SO(3).& &\end{array}
$$
Thus $G_\mu^0$ is a amalgamated free product of two compact subgroups,
$\SO(3)\times \SO(3)$, which is the automorphism group of the product almost 
complex structure, and $ S^1\times \SO(3)$. The latter appears as the
automorphism group of the other integrable almost complex structure 
with K\"ahler form $\om_\mu$, namely the Hirzebruch structure
on ${\mathcal P}(L_2\oplus \C)$ where the line bundle $L_2\to \CP^1$ has Chern 
number $2$.  As mentioned in~\cite{AG}, this description has 
interesting parallels with the structure of some Kac--Moody groups.

McDuff~\cite{Mcrs} 
proves that the homotopy type of $G_\mu^0$ is constant on all intervals 
$(\ell-1, \ell], \ell> 1$.  However, their full homotopy type for $\mu> 2$ is 
not yet understood, and there are only partial results
when $g>0$.   Apart from this there is rather little known about 
the homotopy type of $\Symp(M)$.  There are some results due to
Pinsonnault~\cite{Pin} and Lalonde--Pinsonnault~\cite{LalP} on the 
one point blow up of $S^2\times S^2$ showing that the homotopy type 
of this group also depends on the symplectic area of the exceptional divisor.  
Also Seidel~\cite{Seid2,SEI} 
has done some very 
interesting work on the 
symplectic mapping class group $\pi_0(\Symp(M))$ for certain 
$4$-manifolds, and on the case $M = \CP^m\times \CP^n$.

\section{The Hamiltonian group}

Now consider the
 Hamiltonian subgroup $\Ham(M)$.  It has many special properties:
 it is the commutator subgroup of $\Symp_0(M)$ and is itself 
 a simple group (Banyaga).  
 It also supports a biinvariant metric, the 
 Hofer metric, which gives rise to an interesting geometry. 
 Its elements also have remarkable dynamical
  properties.  
  For example,  according to  Arnold's conjecture 
   (finally proven by  Fukaya--Ono and Liu--Tian based on work by Floer and 
   Hofer--Salamon)
  the number of fixed points of $\phi\in\Ham$ may be estimated as
  $$
  \#{\rm Fix\,} \phi\ge \sum_k {\rm rank\,} H^k(M, \Q)
  $$
  provided that  the fixed points are all nondegenerate, i.e. that the graph  
  of $\phi$ is transverse to the diagonal.  
  
  Many features of this group are still not understood, and it 
may not even be $C^1$-closed in $\Symp_0$. Nevertheless, we will see
that there are some analogs of the stability properties discussed 
earlier for $\Symp$.
Also the action of $\Ham(M)$ on $M$ has special properties.\MS

\subsection*{Hofer Geometry}

Because the elements of the Hamiltonian group are generated by 
functions $H_t$, the group itself supports a variety of interesting functions.
First of all there is the Hofer norm~\cite{Hof} that is usually
defined as follows:
$$
\|\phi\| : = \inf_{\phi_1^H = \phi}\;\int_0^1 \Bigl(\max_{x\in M} H_t(x) - \min_{x\in M} 
H_t(x)\Bigr)dt.
$$
Since this is constant on conjugacy classes and symmetric (i.e. 
$\|\phi\| = \|\phi^{-1}\|$), it gives rise to a 
biinvariant metric $d(\phi,\psi): = \|\psi\phi^{-1}\|$ 
on $\Ham(M,\om)$.   There are still many open 
questions about this norm --- for example, it is not yet known whether 
it is always unbounded: for a good introduction 
see Polterovich's lovely book~\cite{Pbk}.
 
Recently, tools (based on Floer homology) have been developed that allow 
one to define functions on $\Ham$ or its universal cover $\THam$
by picking out special elements of the action spectrum $\Spec(\tphi)$
of $\tphi\in \THam$.  This spectrum is defined as follows.  
Choose a normalized 
time periodic Hamiltonian $H_t$ that generates $\tphi$, i.e. so that 
the following conditions are satisfied:
$$
\int H_t \om^n = 0,\;t\in \R,\qquad H_{t+1}= H_t, \;t\in \R,\qquad 
\tphi = \tphi^H: = (\phi_1^H,\{\phi_t^H\}_{t\in [0,1]}).
$$
Denote by $\TLl(M)$ the cover of the space $\Ll(M)$ of contractible
loops $x$ in $M$ whose elements are pairs $(x, u)$, where $u:D^2\to M$ 
restricts to $x$ on $\p D^2 = S^1$.  Then define the 
action functional  $\Aa_H:\TLl(M) \to \R$ by setting
$$
\Aa_H(x, u) = \int_0^1 H_t(x_t)\,dt - \int_{D^2} u^*(\om).
$$
The critical points of $\Aa_H$ are precisely the pairs $(x,u)$ where 
$x$ is a contractible $1$-periodic orbit of the flow $\phi_t^H$.  
Somewhat surprisingly, it turns out that the set of critical values 
of $\Aa_H$ depends only on the element $\tphi^H\in \THam$ defined by
the flow $\{\phi_t^H\}_{t\in [0,1]}$; in other words, these values 
depend only on the homotopy class of the path $\phi_t^H$ rel endpoints.
Thus we set:
$$
\Spec(\tphi^H): = \{\mbox{all critical values of }\Aa_H \}.
$$
There are variants of the Hofer norm that pick out certain special
homologically visible elements from this spectrum: see for example 
Schwarz~\cite{Sch} and Oh~\cite{Oh2}.  

Even more interesting is a 
recent construction by Entov--Polterovich~\cite{EnP} that uses these 
spectral invariants to define a 
nontrivial continuous and homogeneous {\bf quasimorphism} $\mu$  on 
$\THam(M,\om)$, when $M$ is a monotone manifold such as 
$\CP^n$ that has semisimple  quantum homology ring.  
A quasimorphism on a 
group $\G$ is a map 
$\mu: \G\to \R$  that is a bounded distance away from being a 
homomorphism, i.e. there is a constant $c = c(\mu)> 0$ such that
$$
|\mu(gh) - \mu(g)-\mu(h)| < c, \qquad g,h\in \G.
$$
It is called homogeneous if $\mu(g^m) = m\mu(g)$ for all $m\in \Z$,
in which case it restricts  to a homomorphism on all abelian subgroups.
Besides giving information about the bounded cohomology of $\G$, 
quasimorphisms
can be used to investigate the commutator lengths  and
dynamical properties of its elements.  The example constructed by  
Entov--Polterovich extends the Calabi homomorphism defined on
the subgroups $\THam_U$ of elements with support in sufficiently small 
open sets $U$.  Moreover, in the case of
$\CP^n$, it vanishes on $\pi_1(\Ham)$ and so 
descends to the Hamiltonian group $\Ham$ (which incidentally equals $\Symp_0$ 
since $H^1(\CP^n) = 0$.)     
It is not yet known whether $\THam(M)$ or $\Ham(M)$ 
supports a nontrivial quasimorphism for every $M$.  
Note that these groups have no 
nontrivial homomorphisms to $\R$ because they are perfect.

\subsection*{Relation between $\Ham$ and $\Symp_0$.}

The relation between $\Ham$ and $\Symp_0$ is best understood via the {\bf Flux 
homomorphism.}
Let $\TSymp_0(M)$ denote the universal cover of $\Symp_0(M)$.
Its elements $\tphi$ are equivalence classes of paths $\{\phi_t\}_{t\in [0,1]}$
starting at the identity,
where $\{\phi_t\}\sim \{\phi_t'\}$ iff $\phi_1 = \phi_1'$ and the paths are 
 homotopic rel endpoints.  We define
$$
\Flux (\tphi) = \int_0^1 [\om(\dot{\phi}_t,\cdot)] \in H^1(M, \R).
$$
That this  depends only on the  homotopy class 
of the path $\phi_t$ (rel endpoints) is a consequence of the 
following alternative description:
 the value of the cohomology class $\Flux (\tphi)$ on a 
$1$-cycle
$\ga: S^1\to M$ is given by the integral
\begin{equation}\label{eq:flux0}
\Flux (\tphi)(\ga) = \int_{\tphi_*(\ga)}\om, 
\end{equation}
where $\tphi_*(\ga)$ is the $2$-chain $I\times S^1\to M: 
(t,s)\mapsto \phi_{t}(\ga(s))$.   Thus $\Flux$ is well defined.  It is 
not hard to check that it is a homomorphism.

One of the first results in the theory is that the rows and columns
in the following commutative diagram are short exact sequences of groups.
(For a proof see~\cite[Chapter~10]{MS}.)
\begin{equation}\label{eq:flux}
\begin{array}{ccccc}
 \pi_1(\Ham(M))& \longrightarrow&
\pi_1(\Symp_0(M)) & \stackrel{\Flux}\longrightarrow &  \Ga_\om \\
 \downarrow& &\downarrow & & \downarrow  \\
 \THam(M)& \longrightarrow&
\TSymp_0(M) & \stackrel{\Flux}\longrightarrow&  H^1(M, \R)\\
\downarrow& &\downarrow & & \downarrow\\
\Ham(M)& \longrightarrow & \Symp_0(M) & 
\stackrel{\Flux}\longrightarrow& H^1(M, \R)/\Ga_\om.
\end{array}
\end{equation}
Here $\Ga_\om$ is the so-called {\bf flux group}.  It is the image of
$\pi_1(\Symp_0(M))$ under the flux homomorphism. 

It is easy to see that  $\Ham(M)$ is $C^1$-closed in $\Symp_0(M)$  if 
and only if  $\Ga_\om$ is  a discrete subgroup of $H^1(M, \R)$. 

\begin{question}  Is the subgroup $\Ga_\om$ of $H^1(M, \R)$ always 
discrete?
\end{question}

The hypothesis that $\Ga_{\om}$ is always discrete is known 
as the {\bf Flux conjecture.}  One might think it would always hold 
by analogy with symplectic rigidity.   
In fact it does hold in many cases, for example if $(M,\om)$ is 
K\"ahler or from~(\ref{eq:flux0}) above if $[\om]$ is integral, but
we do not yet have a complete understanding of this question. 
One consequence of Corollary~\ref{cor:rigid}
is that the rank of $\Ga_{\om}$ is always bounded above by the first Betti 
number (see~Lalonde--McDuff--Polterovich~\cite{LMP1,LMP2}; some 
sharper bounds are found in 
Kedra~\cite{Ked1}), but the argument does not rule out the 
possibility that $\Ga_{\om}$ is indiscrete for certain values of 
$[\om]$.  Thus, for the present
 one should think of $\Ham(M)$ as a 
leaf in a foliation of $\Symp_0(M)$ that has 
codimension equal to the first 
Betti number of $M$.

\subsection*{Hamiltonian stability}

When $\Ga_{\om}$ is discrete, the stability principle
extends:   there is a 
$C^1$-neighbourhood of $\Ham(M,\om)$ in $\Diff(M)$ 
that deformation retracts into $\Ham(M,\om)$.  Moreover, if this 
discreteness were uniform with respect to $\om$ (which would hold 
if $(M, \om)$ were K\"ahler),
then the groups 
$\Ham(M, \om)$ would have the same stability with respect to variations in $\om$
as do $\Symp_0$ and $\Symp$.  

To be more precise, suppose that  for
each $\om$ and each $\eps> 0$ there is   a neighbourhood $\Nn(\om)$ such that  
when $\om'\in \Nn(\om)$ $\Ga_{\om'}$ contains no nonzero element 
of norm $\le \eps$.  Then for any compact subset $K$ of $\Ham(M,\om)$
there would be a neighbourhood $\Nn(\om)$ such that $K$ isotops into 
$\Ham(M, \om')$ for each $\om'\in \Nn(\om)$.  For example if $K = 
\{\phi_t\}$ is a 
loop (image of a circle) in $\Ham(M,\om)$
and $\om_s, 0\le s\le 1,$ is any path, this would mean that
any smooth extension 
$\{\phi_t^s\}, s\ge 0,$ of $\{\phi_t\}$ to a family of loops in $\Symp(M, 
\om_s)$ would be homotopic through 
$\om_s$-symplectic loops  to a loop in  $\Ham(M, \om_s)$. 

Even if this hypothesis on $\Ga_\om$ held, it would not rule out 
the possibility of global instability:
a loop in $\Ham(M,\om)$ could be 
isotopic through (nonsymplectic) loops in  $\Diff(M)$ to 
a nonHamiltonian  loop in 
some other far away symplectomorphism group $\Symp(M, \om')$. 
One of the main results in 
Lalonde--McDuff--Polterovich~\cite{LMP2} is that this  global instability
never occurs;
any $\om'$-symplectic loop that is  isotopic in $\Diff(M)$ 
to an $\om$-Hamiltonian loop
must be homotopic in $\Symp(M,\om')$ to an
$\om'$-Hamiltonian  loop regardless of the relation between $\om$ 
and $\om'$ and no matter whether any of the groups  $\Ga_{\om}$ are discrete.  
This is known as {\bf Hamiltonian rigidity} and is a 
consequence of a vanishing theorem for the Flux homomorphism: see 
Corollary~\ref{cor:rigid} below.  As we 
now explain this extends 
to general results about  the 
action of $\Ham(M)$ on $M$.

\subsection*{Action of $\Ham(M)$ on $M$}\label{ss:act}

There are some suggestive but still incomplete results about the action 
of $\Ham(M)$ on $M$. The first result below is folklore.  It is
a consequence of
the proof of the Arnold conjecture, but as we show below (see 
Lemma~\ref{le:sect}) also follows
from a geometric argument.  The second part is due to
 Lalonde--McDuff~\cite{LMh}.   Although the statements are 
 topological in nature, both proofs are based on the existence 
 of the Seidel representation, a deep fact that uses the properties of 
 $J$-holomorphic curves.

\begin{prop}\label{prop:lm0}
\begin{description}\item[(i)]
The evaluation map $\pi_1(\Ham(M)) \to \pi_1(M)$ is zero.
\item[(ii)] The natural action of $H_*(\Ham(M),\Q)$ on $H_*(M,\Q)$ is 
trivial.
\end{description}
\end{prop}

Here the action $\tr_\phi: H_*(M)\to H_{*+k}(M)$
of an element $\phi\in H_k(\Ham(M))$ is defined as follows:

\begin{quote}{\it 
if $\phi$ is represented by the cycle $t\mapsto \phi_t$ for $t\in V^k$
and $c\in H_*(M)$ is represented by $x\mapsto c(x)$ for $x\in C$
then $\tr_\phi(c)$ is represented by the cycle}
$$
V^k\times C \to M: (t,x)\mapsto \phi_t(c(x)).
$$
\end{quote}

\NI
It is just the action on homology induced by the map $\Ham(M)\times 
M\to M$.  It extends to the group
$(M^M)_{id}$  of 
self-maps of $M$ that are homotopic to the identity, and hence 
depends only on the image of $\phi$ in 
$H_k(M^M)_{id}$.  To say it is trivial means that
$$
\tr_\phi(c) = 0\qquad\mbox{whenever } c\in H_i(M), \;i > 0.
$$
Note that this does {\it not} hold for the action of $H_1(\Symp_0(M))$.  
Indeed by (\ref{eq:flux0}) 
the image under the Flux homomorphism of a loop $\la\in \pi_1(\Symp_0(M))$
is simply
\begin{equation}\label{eq:tr}
\Flux(\la)(\ga) = \langle \om, \tr_\la(\ga)\rangle.
\end{equation}

The  rigidity of Hamiltonian loops is an immediate consequence
of Proposition~\ref{prop:lm0}.

\begin{cor}\label{cor:rigid} Suppose that $\phi\in \pi_1(\Symp(M, \om))$ and 
$\phi'\in \pi_1(\Symp(M, \om'))$ represent the same element of 
$\pi_1((M^M)_{id})$.  Then 
$$
\Flux_{\om}(\phi) = 0 \quad\IFF\quad \Flux_{\om'}(\phi') = 0.
$$
\end{cor}

\NI
{\it Proof.}  If $\Flux_{\om}(\phi) = 0$ then $\phi$ is an 
$\om$-Hamiltonian loop and
Proposition~\ref{prop:lm0}(ii) implies that $\tr_\phi: H_1(M)\to 
H_2(M)$
is the zero map. 
But, for each $\ga\in H_1(M)$,~(\ref{eq:tr}) implies that
$$
\Flux_{\om'}(\phi')(\ga) = \Flux_{\om'}(\phi)(\ga) = 
\langle \om', \tr_\phi(\ga)\rangle = 0.
$$
\MS

This corollary is  elementary when the loops are 
 circle subgroups since then one can distinguish between Hamiltonian 
 and nonHamiltonian loops by looking at the weights of the action at 
 the fixed  points: a circle action is Hamiltonian if and only if 
 there is a point whose weights all have the same sign.   One can also 
consider  maps $K\to \Ham(M, \om)$ with arbitrary compact domain $K$. 
But their stability follows from the above result because
  $\pi_k(\Ham(M)) = \pi_k(\Symp_0(M))$ when $k>1$ by
diagram~(\ref{eq:flux}). For more details see~\cite{LMs}.

Thus one can 
compare the homotopy types of the groups $\Ham(M, \om)$ (or of 
$\Symp(M, \om)$) as $[\om]$ 
varies in $H^2(M, \R)$.  More precisely, as Buse points out in~\cite{Bu},
any element $\al$ in $\pi_*(\Ham(M, \om))$ has a smooth extension to a 
family $\al_t\in \pi_*(\Ham(M, \om_t))$ where  $[\om_t]$ fills out a
neighborhood of $[\om_0] = [\om]$ in $H^2(M, \R)$.
Moreover the germ of this extension at $\om = \om_0$ is unique.     
Thus one can distinguish between {\bf robust} elements in the 
homology or homotopy of the spaces $\Ham(M, \om_0)$ and $B
\Ham(M, \om_0)$  whose extensions are nonzero for all $t$ near $0$ 
and {\bf fragile} elements whose extensions vanish as $[\om_t]$ 
moves in certain 
directions.   For example, any class in $H^*(B\Ham(M, \om_0))$ 
that is detected by Gromov--Witten invariants
(i.e. does not vanish on a suitable space of $J$-holomorphic curves as 
in Le--Ono~\cite{LeO}) 
is robust, while the classes $w_\ell$ of 
Proposition~\ref{prop:am2} are fragile.  For some interesting examples 
in this connection, see Kronheimer~\cite{Kro} and Buse~\cite{Bu}.

\subsection*{$c$-splitting for Hamiltonian bundles}

From now on, we assume that (co)homology has rational coefficients.
Since the rational cohomology $H^*(G)$ of any $H$-space  (or group) is 
freely generated by the dual of its rational homotopy, it is easy to
see that part (ii) of Proposition~\ref{prop:lm0} holds if and only if 
it holds for all spherical classes
$\phi\in H_k(\Ham(M))$.  Each such $\phi$ gives rise to a 
locally trivial fiber bundle  $M\to P_{\phi}\to S^{k+1}$ with structural group 
$\Ham(M)$.
Moreover, 
the differential in the corresponding Wang sequence is precisely  
$\tr_\phi$.
In other words, there is an exact sequence:
\begin{equation}\label{eq:wang}
\dots H_{i}(M)\stackrel{\tr_\phi}\to H_{i+k}(M) \to  H_{i+k}(P_{\phi})
\stackrel{\cap [M]}\to
H_{i-1}(M) \to\dots
\end{equation}
Hence $\tr_\phi = 0$ for $k>0$ if and only if this long exact sequence 
breaks up into short exact sequences:
$$
0 \to H_{i+k}(M) \to  H_{i+k}(P_{\phi}) \stackrel{\cap [M]}\to
H_{i-1}(M) \to 0.
$$
Thus Proposition~\ref{prop:lm0}(ii)  is equivalent to the following 
statement.

\begin{prop}\label{prop:ham1} For every Hamiltonian bundle $P\to S^{k+1}$,
with fiber $(M,\om)$ the rational homology $H_*(P)$ 
is isomorphic as a 
vector space to the tensor product $H_*(M)\otimes H_*(S^{k+1})$.  
\end{prop}
 
Observe that the corresponding isomorphism in cohomology
 need not preserve the ring structure.  We say that 
a bundle  $M\to P \to B$  is {\bf c-split} if the rational cohomology
$H^*(P)$ is isomorphic  as a vector space to $H^*(M)\otimes H^*(B)$.

\begin{question} \label{q:csp} Is every 
 fiber bundle $M\to P \to B$ with structural group $\Ham(M)$ c-split?
\end{question}

It is shown in~\cite{LMh} that the answer is affirmative if $B$ has 
dimension $\le 3$ or is a product of spheres and projective spaces 
with fundamental group of rank $\le 3$.  
By an old result of Blanchard, it is also affirmative if 
$(M,\om)$ satisfies the hard Lefschetz 
condition, i.e. if
$$
\wedge [\om]^{k}: \quad H^{n-k}(M,\R) \to
H^{n+k}(M,\R)
$$
is an isomorphism for all $0< k< n$.  (This argument has now been 
somewhat extended by Haller~\cite{Ha} using ideas of Mathieu about 
the harmonic cohomology of a symplectic manifold.)  If the 
structural group of $P\to B$ reduces to a finite dimensional Lie group 
$G$, then $c$-splitting is equivalent to
a result of Atiyah--Bott~\cite{AB} about the 
structure of the equivariant cohomology ring $H_G^*(M)$.  
This is the cohomology of  
 the universal Hamiltonian $G$-bundle with fiber $M$
$$
M \;\longrightarrow\; M_G = EG\times_G M \;\longrightarrow\; BG,
$$
 and was shown in~\cite{AB} to be isomorphic to 
$H^*(M)\otimes H^*(BG)$ as a $H^*(BG)$-module.  Hence a positive
answer to 
Question~\ref{q:csp} in general would imply that this aspect of the homotopy theory of 
Hamiltonian actions is similar to the more rigid cases, when the group is 
finite dimensional or when the manifold is K\"ahler. 
 For further discussion 
see~\cite{LMh,LMs} and Kedra~\cite{Ked2}.

Note finally that all results on the action of $\Ham(M)$ on $M$ can be 
phrased in terms of the universal Hamiltonian bundle 
$$
M \to M_{\Ham} = E\Ham \times_{\Ham} M \to B\Ham(M).
$$
For example, Proposition~\ref{prop:lm0} part (i)  states that this bundle 
has a section over its $2$-skeleton.  Such a formulation has the 
advantage that it immediately suggests further questions.  For 
example, one might wonder if the bundle $M_{\Ham} \to B\Ham$ always 
has a global section.  However this fails when $M = S^2$ since the 
map $\pi_3(\Ham(S^2)) = \pi_3(\SO(3)) \to \pi_3(S^2)$ is nonzero.

\section{Symplectic geometry of bundles over $S^2$}

The proofs of Propositions~\ref{prop:lm0} and \ref{prop:ham1} above rely on 
properties of Hamiltonian bundles over $S^2$.   
We now show how  the Seidel representation
$$
\pi_1(\Ham(M,\om))\to (\QH_{ev}(M))^{\times}
$$
of $\pi_1(\Ham(M,\om))$ into the group of even units in quantum 
homology  gives information
on the homotopy properties of Hamiltonian bundles.
As preparation, we first discuss  quantum homology.

\subsection*{The small quantum homology ring $QH_*(M)$}

There are several slightly different ways of defining 
the small quantum homology ring.  We adopt the conventions of~\cite{LMP2,Mcq}.

Set $c_1 = c_1(TM)\in H^2(M,\Z)$.  Let $\La$ be the Novikov 
ring of
the group $\Hh = H_2^S(M,\R)/\!\!\sim$  with valuation $I_\om$ where 
$B\sim B'$
if $\om(B-B') = c_1(B-B') = 0$.
 Then $\La$ is the completion of the rational group ring
 of $\Hh$ with elements of the form 
$$
\sum_{B\in \Hh} q_B\; e^B
$$
where for each $\ka$ there are only finitely many nonzero
$q_B\in \Q$ with $\om(B) > - \ka$.
  Set 
$$
QH_{*}(M) = QH_*(M,\La)
= H_*(M)\otimes\La.
$$
We may define an $\R$ grading on $QH_*(M,\La)$ by setting
$$
\deg(a\otimes e^B) = \deg(a) + 2c_1(B),
$$
and can also think of
$QH_*(M,\La)$ as $\Z/2\Z$-graded with 
$$
QH_{\ev} =  
H_{\ev}(M)\otimes\La, \quad QH_{odd} =  
H_{odd}(M)\otimes\La.
$$

Recall that the quantum intersection product 
$$
a*b\in QH_{i+j - 2n}(M), \qquad \mbox {for }\; a\in H_i(M), b\in 
H_j(M)
$$
 is defined as follows:
\begin{equation}\label{eq:qm0}
a* b = \sum_{B\in \Hh} (a*b)_B\otimes e^{-B},
\end{equation} 
where  $(a*b)_B\in H_{i+j- 2n+2c_1(B)}(M)$ is defined by the 
requirement
that 
\begin{equation}\label{eq:qm}
(a*b)_B\,\cdot\, c = \GW_M(a,b,c;B) \quad\mbox{ for all }\;c\in 
H_*(M).
\end{equation}
Here  $\GW_{M}(a,b,c; B)$ denotes the Gromov--Witten invariant that counts 
the number of $B$-spheres in $M$ meeting the cycles $a,b,c\in 
H_{*}(M)$,
and we have written $\cdot $ for the usual  intersection
pairing on $H_*(M) = H_*(M, \Q)$.  Thus
$a\cdot  b = 0$ unless $\dim(a) + \dim (b) = 2n$ in which case it is 
the
algebraic number of intersection points of the cycles. 

Alternatively, one can define $a*b$ as follows: if
$\{e_i\}$ is a basis for $H_*(M)$ with dual basis $\{e_i^*\}$, then
$$
a * b = \sum_i \GW_{M}(a,b,e_i; B)\; e_i^*\otimes e^{-B}.
$$
The product  
$*$ is
extended to $QH_*(M)$ by linearity over $\La$, and is associative. 
Moreover,
it  preserves the $\R$-grading in the homological sense, i.e. it obeys 
the same grading rules as does the intersection product.  

 This product $*$ gives $QH_{*}(M)$
the structure of a
graded commutative ring with unit $\1 = [M]$. Further, the invertible
elements in $QH_{\ev}(M)$ form a commutative group  
$(QH_{\ev}(M,\La))^\times$
that acts on   $QH_*(M)$ by  quantum multiplication.  By Poincar\'e 
duality one can transfer this product to cohomology.  Although this is 
very frequently done, it is often easier to work with homology when 
one wants to understand the relation to geometry.

\subsection*{The Seidel representation $\Psi$}

Consider a smooth bundle $\pi: P  \to 
S^2$ with fiber $M$.
Here we consider $S^2$ to be the union $D_+ \cup D_-$ of two copies of 
$D$, with the  orientation  of $D_+$.  We denote the equator $ D_+\cap 
D_-$ by $\p$, oriented as the boundary of $D_+$, and choose some point
$*$ on $\p$ as the base point of $S^2$.  We assume also that 
the fiber $M_*$ over $*$ has a chosen identification with $M$.

Since every smooth bundle over a disc can be trivialized, 
we can build any smooth bundle $P\to S^2$ by taking two product 
manifolds $D_\pm\times M$ and gluing them along the boundary $\p\times M$ by a 
based 
loop $\la = \{\la_t\}$ in $\Diff(M)$.  Thus
$$
P = (D_+\times M)\; \cup \;(D_-\times M)/\sim,\quad (e^{2\pi it},x)_- \equiv 
(e^{2\pi it},\la_t(x))_+.
$$
A {\bf symplectic bundle} is built from a based loop in $\Symp(M)$ and a 
{\bf Hamiltonian bundle} from one in $\Ham(M)$.  Thus the smooth bundle
$P\to S^2$ is symplectic if and only if there is a smooth family of 
cohomologous symplectic 
forms $\om_b$ on the fibers $M_b$.   
It is shown in~\cite{Seid,MS,LMh} 
that a symplectic bundle $P\to S^2$ is Hamiltonian if and only if
the fiberwise forms $\om_b$ have a closed extension $\Om$.  
(Such forms $\Om$ are called {\bf $\om$-compatible}.)
 Note that in any of these categories
two bundles are equivalent if and only if their defining loops are 
homotopic.

From now on, we restrict to Hamiltonian bundles, and denote by
$P_\la\to S^2$   the bundle constructed from a loop $\la\in 
\pi_1(\Ham(M))$.
By adding the pullback of a 
suitable area form on the base we can choose the closed extension
 $\Om$ to be  symplectic. 
The manifold $P_\la$ carries two canonical 
cohomology classes,
the first Chern class of the vertical
tangent bundle 
$$ 
c_{vert} = c_1(TP_\la^{vert})  \in H^2(P_{\la},\Z),
$$
and the   coupling class  $u_\la$, i.e.  the unique class
in $H^2(P_\la,\R)$ such that 
$$
i^*(u_\la) = [\om],\qquad u_\la^{n+1} = 0,  
$$
where $i: M\to P_\la$ is the inclusion of a fiber.

The next step is to choose a canonical (generalized) section class
in  $\si_{\la}\in H_{2}(P_{\la}, 
 \R)/\sim$.  By definition this should project onto the positive generator
of $H_{2}(S^{2}, \Z)$,
 In the general case, 
when $c_{1}$ and $[\om]$ induce linearly 
independent homomorphisms $H_2^S(M) \to \R$, $\si_\la$ is
defined by the requirement that
\begin{equation}\label{eq:seccl}
c_{vert}(\si_{\la}) \;= u_{\la}(\si_{\la})\; =\;0,
\end{equation}
which has a unique solution modulo the given equivalence.
If either $[\om]$ or $c_1$ vanishes on 
 $H_2^S(M)$ then
such a class $\si_\la$ still exists.\footnote
{
See~\cite[Remark~3.1]{Mcq} for 
the case when $[\om] = 0$ on  $H_2^S(M)$.  If $c_1 = 0$ on  $H_2^S(M)$
but $[\om]\ne 0$ then we can choose
$\si_{\la}$ so that $u_\la(\si_{\la}) = 0$.  Since 
$c_{vert}$ is constant on section classes, we must show that it 
always vanishes.  But the 
existence of the Seidel representation implies
every Hamiltonian fibration $P\to S^2$ has 
some section $\si_P$ with $n\le c_{vert}(\si_P) \le 0$ (since it only 
counts such sections), and the value 
must be $0$ because 
$c_{vert}(\si_{P_\la\#P_{-\la}}) = c_{vert}(\si_{P_\la})+ 
c_{vert}(\si_{P_{-\la}})$:
see~\cite[Lemma~2.2]{Mcq}.}
In the remaining case (the monotone case), when $c_1$ is some nonzero multiple of 
$[\om]\ne 0$ 
on $H_2^S(M)$,
we choose $\si_\la$ so that $c_{vert}(\si_\la) = 0$.

We then set
\begin{equation}\label{eq:qm2}
\Psi(\la) = \sum_{B\in \Hh} a_B\otimes e^{B}
\end{equation}
where, for all $c\in H_{*}(M)$,
\begin{equation}\label{eq:qm3}
a_{B}\cdot_{M} c = \GW_{P_{\la}}([M], [M], c\,; \si_{\la} - B).
\end{equation}
Note that $\Psi(\la)$ belongs to the strictly commutative part 
$QH_{ev}$ of $QH_{*}(M)$.
Moreover $\deg(\Psi(\la)) = 2n$ because $c_{vert}(\si_{\la}) = 0.$
Since all $\om$-compatible forms are deformation 
equivalent, $\Psi$ is independent of the choice of $\Om$.

Here is the  main result.

\begin{prop}\label{prop:seid}
For all $\la_1,\la_2\in \pi_1(\Ham(M))$
$$
\Psi(\la_1 + \la_2) = \Psi(\la_1)*\Psi(\la_2),\qquad 
\Psi(0) = \1,
$$
where $0$ denotes the constant loop. Hence 
$\Psi(\la)$ is invertible for all $\la$ and $\Psi$ defines a 
group homomorphism
$$
 \Psi: \pi_{1}(\Ham(M,\om))\;\;\to\;\; (QH_{\ev}(M,\La))^{\times}.
$$
\end{prop}

In the case when $(M, \om)$ 
satisfies a suitable positivity condition, 
this is  a variant of the main result in 
Seidel~\cite{Seid}.  The general proof is
due to McDuff~\cite{Mcq} 
using ideas from Lalonde--McDuff--Polterovich~\cite{LMP2}. It uses a 
refined version of the ideas in the proof of Lemma~\ref{le:sect} below.  

\subsection*{Homotopy theoretic consequences of the existence of $\Psi$}

First of all, note that because $\Psi(\la)\ne 0$ there must always be 
$J$-holomorphic sections of $P_\la\to S^2$ to count.  Thus every 
Hamiltonian bundle $\pi:P\to S^2$ must have a section $S^2\to 
P$.  If we trivialize $P$ over the two hemispheres $D_{\pm}$ of $S^2$   
and homotop the section to be constant over one of the discs, it becomes 
clear that there is a section  if and only if the defining loop $\la$ of $P$ has 
trivial image under the  evaluation map $\pi_1(\Ham(M))\to \pi_1(M)$.
This proves part (i) of Proposition~\ref{prop:lm0}.

In fact one does not need the full force of Proposition~\ref{prop:seid} 
in order to arrive at this conclusion, since we only have to produce 
one section.  

\begin{lemma}\label{le:sect}
Every Hamiltonian bundle $P\to S^2$ has a section.    
\end{lemma}
\NI
{\it Sketch of Proof} Let $\la = \{\la_t\}$ be a Hamiltonian loop and
consider the family of trivial bundles
$P_{\la, R}\to S^2$ given by
$$
P_{\la, R} = (D_+\times M)\;\cup\; (S^1\times [-R,R]\times M)\;\cup\; 
(D_-\times M)
$$
with attaching maps
$$
\left(e^{2\pi it},\la_t(x)\right)_+ \equiv \left(e^{2\pi it}, -R,x\right),\qquad 
(e^{2\pi it}, R, x)\equiv \left(e^{2\pi it}, \la_t(x)\right)_-.
$$
Thus, $P_{\la, R}$ can be thought of as the fiberwise union (or Gompf 
sum)  of 
$P_{\la}$ with $P_{-\la}$ over a neck of length $R$. It is possible to 
define a family $\Om_R$ of $\om$-compatible symplectic forms on $P_{\la, R}$
in such a way that the manifolds 
$(P_{\la, R}, \Om_R)$ converge in a well defined sense as $R\to 
\infty$.  The limit 
is a singular manifold $P_{\la}\cup P_{-\la}\;\to\; S_{\infty}$ that is 
a locally trivial fiber bundle over the nodal curve 
consisting of the 
one point union of two $2$-spheres.  To do this, one first models the 
convergence of the $2$-spheres in the base by a $1$-parameter family 
$S_R$ of  
disjoint holomorphic spheres in 
the one point blow up of $S^2\times S^2$ that 
converge to the pair $S_{\infty}=\Si_+\cup \Si_-$ of 
exceptional divisors at the blow up point.
Then one builds a suitable smooth Hamiltonian bundle 
$$
\pi_X: (\Xx, 
\Tilde{\Om})\to \Ss
$$
with fiber $(M,\om)$, where  $\Ss$ is a 
neighbourhood  of $\Si_+\cup \Si_-$ in the blow up that contains the union 
of the spheres $S_R, R\ge R_0$:
see~\cite{Mcq}~\S2.3.2.  
The almost complex 
structures $\TJ$ that one puts on $\Xx$ should be chosen so that the 
projection to $\Ss$ is holomorphic.  Then each submanifold $P_{\la, 
R}: = \pi_X^{-1}(S_R)$ 
is $\TJ$-holomorphic.

The bundles $(P_{\la, R}, \Om_R)\to S^2$ are all trivial, and hence
there is one $\TJ$-holomorphic curve in the class $\si_0=[S^2\times pt]$ 
through each point $q_R\in P_{\la, R}$.  
(It is more correct to say that the corresponding 
Gromov--Witten invariant $GW_{P_{\la, R}}([M], [M], pt;\si_0) $ is one; i.e. 
one counts the curves with appropriate multiplicities.) 
Just as in gauge theory, these curves do not disappear when one 
stretches the neck, i.e. lets $R\to \infty$.  Therefore 
as one moves the point $q_R$ to the singular fiber
the family of $\TJ$-holomorphic curves through $q_R$ converges to 
some cusp-curve (stable map)  $C_{\infty}$ in the limit. Moreover, 
 $C_{\infty}$  must lie entirely in the singular fiber 
 $P_{\la}\cup P_{-\la}$ and projects 
to a holomorphic curve in $\Ss$ in the class $[\Si_+]+[\Si_-]$.   Hence 
it must have at least two components, one a section of 
$P_{\la}\to\Si_+$ and the other a section of $P_{-\la}\to\Si_+$.  
There might also be some bubbles in the $M$-fibers, but this is 
irrelevant.\QED

The above argument is relatively easy, in that it
only uses the compactness theorem for 
$J$-holomorphic curves and not the more subtle gluing arguments 
needed to prove things like the associativity of quantum multiplication. 
However the proof of the  rest of  Proposition~\ref{prop:lm0} is based on the fact 
that each element $\Psi(\la)$ is a multiplicative  unit in quantum homology.
The only known way to prove this is via some sort of gluing argument.
 Hence in this case it seems that one does need the full 
force of the gluing arguments, whether one works as here with 
$J$-holomorphic spheres or as in Seidel~\cite{Seid} with Floer homology.

We now show how to deduce part (ii) of Proposition~\ref{prop:lm0}
from Proposition~\ref{prop:seid}.
So far, we have described $\Psi(\la)$ as a unit in $QH_*(M)$.
This unit induces an automorphism of $QH_*(M)$ by quantum 
multiplication on the left:
$$
b\mapsto \Psi(\la) * b,\quad b\in QH_*(M).
$$
The next lemma shows that when $b\in H_*(M)$ then the element 
$\Psi(\la) * b$ can also be described by counting curves in $P_\la$ 
rather than in the fiber $M$.

\begin{lemma} If
$\{e_i\}$ is a basis for $H_*(M)$ with dual basis $\{e_i^*\}$, then
$$
\Psi(\la) * b = \sum_i \GW_{P_\la}([M],b,e_i; \si_\la - B)\; e_i^*\otimes e^B.
$$
\end{lemma}
{\it Sketch of Proof:}
To see this, one first shows that for any section class $\si$
 the invariant $\GW_{P_\la}([M],b,c; \si)$
may be calculated using a fibered $J$  (i.e. one 
for which the projection $\pi:P\to S^2$ is 
holomorphic) and with representing cycles for $b,c$ that each lie in a fiber.  
Then one is counting sections of $P\to S^2$.  If the representing cycles 
for $b,c$ are moved into the same fiber, then the curves must degenerate.
Generically the limiting stable map will have two components, a section in 
some class $\si- C$ together with a $C$ curve that meets $b$ and 
$c$.  Thus, using much the same arguments that prove the usual $4$-point 
decomposition rule, one shows that
\begin{equation}\label{eq:mult}
\GW_{P_\la}([M],b,c; \si) = \sum_{A,i} \GW_{P_\la}([M],[M],e_i; \si - A)
\cdot \GW_M(e_i^*,b,c; A).
\end{equation}
But  $\Psi(\la) = \sum q_j e_j^*\otimes e^B$ where 
$$
q_j = \GW_{P_\la}([M], [M], e_j; \si_\la - B) \in \Q.
$$
Therefore 
\begin{eqnarray*}
\Psi(\la) * b & = & \sum_{C,k} \GW_M(\Psi(\la),b, e_k; C)\; e_k^*\otimes e^{-C} \\
& = & \sum_{B,C,j,k} \GW_{P_\la}([M], [M], e_j; \si_\la - B)\cdot
 \GW_M(e_j^*, b,e_k; C) \;e_k^*\otimes e^{B-C}\\
& = & \sum_{A,k} \GW_{P_\la}([M],b,e_k;\si_\la - A)\; e_k^*\otimes e^A
\end{eqnarray*}
where the first equality uses the definition of $*$, 
the second uses the definition of $\Psi(\la)$
  and the third uses~(\ref{eq:mult}) with $\si = \si_\la - (B-C)$.
  For more details, see~\cite[Prop~1.2]{Mcq}. \QED

Since $\Psi(\la) $ is a unit, the map $b\mapsto \Psi(\la) * b$ is injective.
Hence 
for every $b\in H_*(M)$ there has to be some nonzero invariant
$\GW_{P_\la}([M],b,c; \si_\la - B)$ in $P_\la$.
  In particular, the image $i_*(b)$ of 
the class $b$ in $H_*(P_\la)$ cannot vanish.  Thus the map 
$$
i_*: H_*(M)\to H_*(P_\la)
$$
of rational homology groups is injective.  
By~(\ref{eq:wang}), this implies that 
the  homology of $P_\la$ is isomorphic to the tensor product 
$H_*(S^2)\otimes H_*(M)$.  Equivalently, the map 
$$
\tr_\la: H_*(M)\to H_{*+1}(M)
$$
is identically zero.  This proves Proposition~\ref{prop:lm0} (ii) in the 
case of loops.  The proof for the higher homology $H_*(\Ham)$ with $*>1$
is purely topological.  Since $H^*(\Ham)$ is generated multiplicatively by
elements dual to the homotopy, one first reduces to the case when 
$\phi\in \pi_k(\Ham)$.  Thus we need only see that
 all Hamiltonian bundles $M\to P\to B$ with base $B = S^{k+1}$ are 
 c-split, i.e. that 
Proposition~\ref{prop:ham1} holds.  
Now observe:

\begin{lemma} \begin{description}\item[(i)]  
Let $M\to P'\to B'$ be the pullback of $M\to P\to B$ by a map
$B'\to B$ that induces a surjection on rational homology.  Then 
if $M\to P'\to B'$ is c-split, so is $M\to P\to B$.

\item[(ii)]  Let $F\to X\to B$ be a Hamiltonian bundle in which $B$ is simply 
connected.  Then if all Hamiltonian bundles over $F$ and over $B$ are 
c-split, the same is true for Hamiltonian bundles over $X$.
\end{description}\end{lemma}

(The proof is easy and is given in~\cite{LMh}.) 
This lemma implies that in order to establish c-splitting 
when $B$ is an arbitrary sphere  it suffices 
to consider the cases $B = \C P^n$, $B =$, 
the $1$-point blow up $X_n$ of $\C P^n$,  and $B = T^2\times \C P^n$.
But the first two cases  can be proved by induction using the 
lemma above and the Hamiltonian bundle
$$
\C P^1 \to X_n \to \CP^{n-1},
$$
and the third follows by considering the trivial bundle
$$
T^2 \to T^2\times \C P^n \to  \C P^n.
$$
This completes the proof of Proposition~\ref{prop:ham1}.
Though these arguments can be somewhat extended, they do not seem powerful 
enough to deal with all Hamiltonian bundles.   For some further work 
in this direction, see Kedra~\cite{Ked2}.


\begin{thebibliography}{999999999} 
  
\bibitem{Abr} M. Abreu, Topology of symplectomorphism groups of
    $S^2\times S^2$, {\it Inv. Math.}, {\bf 131} (1998), 1-23.

\bibitem{A} S. Anjos,  The homotopy type of symplectomorphism groups of
    $S^2\times S^2$, {\it Geometry and Topology}, {\bf 6} (2002), 195--218.

\bibitem{AG} S. Anjos and G. Granja, Homotopy decomposition of a  group
    of symplectomorphisms of  $S^2\times S^2$,  AT/0303091.
  
\bibitem{AB}  M. F. Atiyah and R. Bott, The moment map and equivariant
     cohomology. {\it Topology\/}, {\bf 23} (1984), 1--28. 


\bibitem{AM} M.~Abreu and D.~McDuff, Topology of symplectomorphism groups  
    of rational ruled surfaces, SG/9910057, {\it Journ. of Amer. Math. Soc.},
    {\bf 13}, (2000) 971--1009.
  

\bibitem{Bu} O. Buse, Relative family Gromov--Witten invariants and 
     symplectomorphisms, SG/01110313.
 
\bibitem{EnP} M. Entov and L. Polterovich, Calabi quasimorphism and 
quantum homology, SG/0205247, {\it International Mathematics Research 
Notes} (2003).

\bibitem{GRO}    M. Gromov,  Pseudo holomorphic curves in symplectic manifolds, 
       {\it  Inventiones Mathematicae\/}, {\bf 82} (1985), 307--47. 

\bibitem{Ha}  S. Haller, Harmonic cohomology of symplectic manifolds, preprint 
     (2003)

\bibitem{Hof}
        H. Hofer, On the topological properties of symplectic maps. {\it 
        Proceedings of the
        Royal Society of Edinburgh\/}, {\bf 115} (1990), 25--38. 

\bibitem{Ked1} J. Kedra, Remarks on the flux groups, {\it 
       Mathematical Research Letters} (2000).

\bibitem{Ked2} J. Kedra, Restrictions on symplectic fibrations, 
       SG/0203232.

\bibitem{Kro} P. Kronheimer, Some nontrivial families of symplectic structures, 
       preprint (1998).

\bibitem{LMc}  F. Lalonde and D McDuff, The classification of ruled
       symplectic $4$-manifolds, {\it Math. Research Letters} {\bf 3},
       (1996), 769--778.

\bibitem{LMh} F. Lalonde and D McDuff, Symplectic structures on 
       fiber bundles, SG/0010275, {\it Topology} {\bf 42} (2003), 
      309--347.

\bibitem{LMs} F. Lalonde and D McDuff, 
        Cohomological properties of 
        ruled symplectic structures, SG/0010277,  in
        {\it Mirror symmetry and string geometry}, ed Hoker, 
        Phong, Yau,  {\it CRM Proceedings and Lecture Notes}, Amer Math Soc. (2001)

\bibitem{LMP1}  F. Lalonde, D. McDuff and L. Polterovich, On the Flux
         conjectures, {\it CRM Proceedings and Lecture Notes} 
        {\bf 15} Amer Math Soc.,  (1998), 69--85.

\bibitem{LMP2}  F. Lalonde, D. McDuff and L. Polterovich,
         Topological rigidity of Hamiltonian
         loops and quantum homology, {\it Invent. Math} {\bf 135}, 369--385 
        (1999)

\bibitem{LalP}   F. Lalonde and M. Pinsonnault, 
      The topology of the space of symplectic balls in rational 
      $4$-manifolds, preprint 2002,  SG/0207096.

\bibitem{LeO}  H. V. Le and K. Ono, Parameterized Gromov--Witten invariants 
      and topology of symplectomorphism groups,
      preprint \#28, MPIM Leipzig (2001) 


\bibitem{Mcq}  D. McDuff, Quantum homology of Fibrations over $S^2$, 
      {\it International Journal of Mathematics}, {\bf 11}, (2000), 665--721.

\bibitem{Mcrs}  D. McDuff,
       Symplectomorphism Groups and Almost Complex Structures,
       SG/0010274,  {\it Enseignment Math.} {\bf 38} (2001), 1--30.

\bibitem{JHOL}  D. McDuff and D.A. Salamon,  {\it $J$-holomorphic curves 
        and Symplectic Topology}, American 
        Mathematical Society, Providence, RI, to appear.  

\bibitem{MS}  D. McDuff and D. Salamon, {\it Introduction to 
      Symplectic Topology}, 2nd edition (1998) OUP, Oxford, UK

     
\bibitem{Oh2}  Yong-Geun Oh, Minimax theory, spectral invariants and 
      geometry of the Hamiltonian diffeomorphism group, SG/0206092.


\bibitem{Pin} M. Pinsonnault,  Remarques sur la Topologie du Groupe 
      des Automorphismes Symplectiques de l'\'Eclatment de $S^2\times S^2$, 
      Ph. D. Thesis, UQAM Montreal (2001).


\bibitem{Pbk} L. Polterovich, {\it The Geometry of the Group of Symplectic 
       Diffeomorphisms}, Lectures in Math, ETH, Birkhauser (2001)

\bibitem{Sch} M. Schwarz, On the action spectrum for closed 
       symplectially aspherical manifolds,
       {\it Pacific Journ. Math} {\bf 193} (2000), 419--461.


\bibitem{Seid}
       P. Seidel, $\pi_1$  of symplectic automorphism groups
       and invertibles in quantum cohomology rings, {\it Geom. and Funct. 
       Anal.} {\bf 7} (1997), 1046 -1095.

\bibitem{SEI} P. Seidel,  On the group of symplectic automorphisms of 
      $\C P^m\times \C P^n$,  {\it Amer. Math. Soc. Transl.} 
      (2) {\bf 196} (1999), 237--250.
 
\bibitem{Seid2} P. Seidel, Graded Lagrangian submanifolds, {\it Bull Math. Soc. 
       France} {\bf  128} (2000), 103--149

\end{thebibliography}
\end{document}